\magnification\magstep1
\input amssym.def
\input amssym.tex
\newread\AUX\immediate\openin\AUX=\jobname.aux
\def\ref#1{\expandafter\edef\csname#1\endcsname}
\ifeof\AUX\immediate\write16{\jobname.aux gibt es nicht!}\else
\input \jobname.aux
\fi\immediate\closein\AUX
\def\today{\number\day.~\ifcase\month\or
  Januar\or Februar\or M\"arz\or April\or Mai\or Juni\or
  Juli\or August\or September\or Oktober\or November\or Dezember\fi
  \space\number\year}
\font\sevenex=cmex7
%\font\sevenex=cmex10 scaled 700
\scriptfont3=\sevenex
%\font\fiveex=cmex7 scaled 714
\font\fiveex=cmex10 scaled 500
\scriptscriptfont3=\fiveex
\def\epsilon{\varepsilon}
\def\square{{\vcenter{\vbox{\hrule height .4pt
          \hbox{\vrule width .4pt height 5.5pt \kern 5.5pt
                \vrule width .4pt}
                            \hrule height .4pt} } }}
\def\theta{\vartheta}

\def\uauf{\lower1.7pt\hbox to 3pt{%
\vbox{\offinterlineskip
\hbox{\vbox to 8.5pt{\leaders\vrule width0.2pt\vfill}%
\kern-.3pt\hbox{\lams\char"76}\kern-0.3pt%
$\raise1pt\hbox{\lams\char"76}$}}\hfil}}
%%%%%%%%%%%%%%%%%%
% Makros f"ur Querverweise:
\def\cite#1{\expandafter\ifx\csname#1\endcsname\relax
{\bf?}\immediate\write16{#1 is not defined!}\else\csname#1\endcsname\fi}
\def\expandwrite#1#2{\edef\next{\write#1{#2}}\next}
\def\neverexpand{\noexpand\noexpand\noexpand}
\def\strip#1\ {}
\def\ncite#1{\expandafter\ifx\csname#1\endcsname\relax
{\bf?}\immediate\write16{#1 is not defined!}\else
\expandafter\expandafter\expandafter\strip\csname#1\endcsname\fi}
\newwrite\AUX
\immediate\openout\AUX=\jobname.aux
%%%%%%%%%%%%%%%%%%%%%
%Macro for numbering formulas, write it in the aux file
\newcount\formula 
\def\eqn#1{\global\advance\formula by 1
\edef\test{\number\formula}
\expandafter\ifx\csname#1\endcsname\test\relax\else
\immediate\write16{#1 has now been defined!}\fi
\expandwrite\AUX{\neverexpand\ref{#1}{(\test)}}
\eqno{(\test)}
}
 
%%%%%%%%%%%%%%%%%%%%%
\newcount\Abschnitt\Abschnitt0
\def\beginsection#1. #2 \par{\advance\Abschnitt1%
\vskip0pt plus.10\vsize\penalty-250
\vskip0pt plus-.10\vsize\bigskip\vskip\parskip
\edef\TEST{\number\Abschnitt}
\expandafter\ifx\csname#1\endcsname\TEST\relax\else
\immediate\write16{#1 has already been defined!}\fi
\expandwrite\AUX{\neverexpand\ref{#1}{\TEST}}
\leftline{\bf\number\Abschnitt. \ignorespaces#2}%
\nobreak\smallskip\noindent\SATZ1}
%%%%%%%%%%%%%%%%%%
\def\Proof:{\par\noindent{\it Proof:}}
\def\Remark:{\ifdim\lastskip<\medskipamount\removelastskip\medskip\fi
\noindent{\bf Remark:}}
\def\Remarks:{\ifdim\lastskip<\medskipamount\removelastskip\medskip\fi
\noindent{\bf Remarks:}}
\def\Definition:{\ifdim\lastskip<\medskipamount\removelastskip\medskip\fi

\noindent{\bf Definition:}}
\def\Example:{\ifdim\lastskip<\medskipamount\removelastskip\medskip\fi
\noindent{\bf Example:}}
%%%%%%%%%%%%%%%%
\newcount\SATZ\SATZ1
\def\proclaim #1. #2\par{\ifdim\lastskip<\medskipamount\removelastskip
\medskip\fi
\noindent{\bf#1.\ }{\it#2}\Par
\ifdim\lastskip<\medskipamount\removelastskip\goodbreak\medskip\fi}
\def\Aussage#1{%
\expandafter\def\csname#1\endcsname##1.{\ifx?##1?\relax\else
\edef\TEST{#1\penalty10000\ \number\Abschnitt.\number\SATZ}
\expandafter\ifx\csname##1\endcsname\TEST\relax\else
\immediate\write16{##1 had already been defined!}\fi
\expandwrite\AUX{\neverexpand\ref{##1}{\TEST}}\fi
\proclaim {\number\Abschnitt.\number\SATZ. #1\global\advance\SATZ1}.}}
\Aussage{Theorem}
\Aussage{Proposition}
\Aussage{Corollary}
\Aussage{Conjecture}
\Aussage{Lemma}
%%%%%%%%%%%%%%%%
\font\la=lasy10
\def\strich{\hbox{$\vcenter{\hbox
to 1pt{\leaders\hrule height -0,2pt depth 0,6pt\hfil}}$}}
\def\dashedrightarrow{\hbox{%
\hbox to 0,5cm{\leaders\hbox to 2pt{\hfil\strich\hfil}\hfil}%
\kern-2pt\hbox{\la\char\string"29}}}

\def\Bindestrich{\penalty10000-\hskip0pt}
\let\_=\Bindestrich
\def\.{{\sfcode`.=1000.}}
%%%%%%%%%%%%%%%%%%%%%%%%%%%%%%%%%%%%

\def\Par{\par}
\def\:={\mathrel{\raise0,9pt\hbox{.}\kern-2,77779pt
\raise3pt\hbox{.}\kern-2,5pt=}}
\def\=:{\mathrel{=\kern-2,5pt\raise0,9pt\hbox{.}\kern-2,77779pt
\raise3pt\hbox{.}}}

\def\|#1|{\mathop{\rm#1}\nolimits}
\def\<{\langle}
\def\>{\rangle}
\let\Times=\times
\def\times{\mathop{\Times}}
\let\Otimes=\otimes
\def\otimes{\mathop{\Otimes}}
%%%%%%%%%%%%%%%%%%%%%%%%%%%%%%%%%
%Laden von Fonts:
\catcode`\@=11
\def\hex#1{\ifcase#1 0\or1\or2\or3\or4\or5\or6\or7\or8\or9\or A\or B\or
C\or D\or E\or F\else\message{Warnung: Setze hex#1=0}0\fi}
\def\fontdef#1:#2,#3,#4.{%
\alloc@8\fam\chardef\sixt@@n\FAM
\ifx!#2!\else\expandafter\font\csname text#1\endcsname=#2
\textfont\the\FAM=\csname text#1\endcsname\fi
\ifx!#3!\else\expandafter\font\csname script#1\endcsname=#3
\scriptfont\the\FAM=\csname script#1\endcsname\fi
\ifx!#4!\else\expandafter\font\csname scriptscript#1\endcsname=#4
\scriptscriptfont\the\FAM=\csname scriptscript#1\endcsname\fi
\expandafter\edef\csname #1\endcsname{\fam\the\FAM\csname text#1\endcsname}
\expandafter\edef\csname hex#1fam\endcsname{\hex\FAM}}
\catcode`\@=12 

%%%%%%%%%%%%%%%%%%%%%%%%%%%%%%%%%
\fontdef Ss:cmss10,,.
\fontdef Fr:eufm10,eufm7,eufm5.

\def\fg{{\Fr g}}
\def\fh{{\Fr h}}
			%Hier aufpassen!!!

\def\fl{{\Fr l}}

\def\fs{{\Fr s}}

\newread\AUXX
\immediate\openin\AUXX=msxym.tex
\ifeof\AUXX
\fontdef bbb:msbm10,msbm7,msbm5.
\fontdef mbf:cmmib10,cmmib7,.
\else
\fontdef bbb:msym10,msym7,msym5.
\fontdef mbf:cmmib10,cmmib10 scaled 700,.
\fi
\immediate\closein\AUXX

\def\FF{{\bbb F}}

\def\RR{{\bbb R}}

\def\ZZ{{\bbb Z}}

\def\cO{{\cal O}}\def\cP{{\cal P}}
\def\cR{{\cal R}}

\mathchardef\leer=\string"0\hexbbbfam3F
\mathchardef\subsetneq=\string"3\hexbbbfam24
\mathchardef\semidir=\string"2\hexbbbfam6E
\mathchardef\dirsemi=\string"2\hexbbbfam6F
\let\OL=\overline
\def\overline#1{{\hskip1pt\OL{\hskip-1pt#1\hskip-1pt}\hskip1pt}}

%<--                    Aufpassen  

%
%%%%%%%%%%%%
% Displayroutine
\abovedisplayskip 9.0pt plus 3.0pt minus 3.0pt
\belowdisplayskip 9.0pt plus 3.0pt minus 3.0pt
\newdimen\Grenze\Grenze2\parindent\advance\Grenze1em
\newdimen\Breite
\newbox\DpBox
\def\NewDisplay#1$${\Breite\hsize\advance\Breite-\hangindent
\setbox\DpBox=\hbox{\hskip2\parindent$\displaystyle{#1}$}%
\ifnum\predisplaysize<\Grenze\abovedisplayskip\abovedisplayshortskip
\belowdisplayskip\belowdisplayshortskip\fi
\global\futurelet\nexttok\WEITER}
\def\WEITER{\ifx\nexttok\qed\expandafter\leftQEDdisplay
\else\leftdisplay\fi}
\def\WEITER\leftdisplay
\def\leftdisplay{\hskip-\hangindent\leftline{\box\DpBox}$$}
\def\leftQEDdisplay{\hskip-\hangindent
\line{\copy\DpBox\hfill\lower\dp\DpBox\copy\QEDbox}%
\belowdisplayskip0pt$$\bigskip\let\nexttok=}
%\everydisplay{\NewDisplay}
%%%%%%%%%%%%
\newbox\QEDbox
\newbox\nichts\setbox\nichts=\vbox{}\wd\nichts=2mm\ht\nichts=2mm
\setbox\QEDbox=\hbox{\vrule\vbox{\hrule\copy\nichts\hrule}\vrule}
\def\qed{\leavevmode\unskip\hfil\null\nobreak\hfill\copy\QEDbox\medbreak}
%%%%%%%%%%%%%%
\newdimen\HIindent
\newbox\HIbox
\def\setHI#1{\setbox\HIbox=\hbox{#1}\HIindent=\wd\HIbox}
\def\HI#1{\par\hangindent\HIindent\hangafter=0\noindent\leavevmode
\llap{\hbox to\HIindent{#1\hfil}}\ignorespaces}
%%%%%%%%%%%%%%

\baselineskip12pt
\parskip2.5pt plus 1pt
\hyphenation{Hei-del-berg}
\def\L|Abk:#1|Sig:#2|Au:#3|Tit:#4|Zs:#5|Bd:#6|S:#7|J:#8||{%
\edef\TEST{[#2]}
\expandafter\ifx\csname#1\endcsname\TEST\relax\else
\immediate\write16{#1 hat sich geaendert!}\fi
\expandwrite\AUX{\neverexpand\ref{#1}{\TEST}}
\HI{[#2]}
\ifx-#3\relax\else{#3}: \fi
\ifx-#4\relax\else{#4}{\sfcode`.=3000.} \fi
\ifx-#5\relax\else{\it #5\/} \fi
\ifx-#6\relax\else{\bf #6} \fi
\ifx-#8\relax\else({#8})\fi
\ifx-#7\relax\else, {#7}\fi\Par}

\def\B|Abk:#1|Sig:#2|Au:#3|Tit:#4|Reihe:#5|Verlag:#6|Ort:#7|J:#8||{%
\edef\TEST{[#2]}
\expandafter\ifx\csname#1\endcsname\TEST\relax\else
\immediate\write16{#1 hat sich geaendert!}\fi
\expandwrite\AUX{\neverexpand\ref{#1}{\TEST}}
\HI{[#2]}
\ifx-#3\relax\else{#3}: \fi
\ifx-#4\relax\else{#4}{\sfcode`.=3000.} \fi
\ifx-#5\relax\else{(#5)} \fi
\ifx-#7\relax\else{#7:} \fi
\ifx-#6\relax\else{#6}\fi
\ifx-#8\relax\else{ #8}\fi\Par}

\def\Pr|Abk:#1|Sig:#2|Au:#3|Artikel:#4|Titel:#5|Hgr:#6|Reihe:{%
\edef\TEST{[#2]}
\expandafter\ifx\csname#1\endcsname\TEST\relax\else
\immediate\write16{#1 hat sich geaendert!}\fi
\expandwrite\AUX{\neverexpand\ref{#1}{\TEST}}
\HI{[#2]}
\ifx-#3\relax\else{#3}: \fi
\ifx-#4\relax\else{#4}{\sfcode`.=3000.} \fi
\ifx-#5\relax\else{In: \it #5}. \fi
\ifx-#6\relax\else{(#6)} \fi\PrII}
\def\PrII#1|Bd:#2|Verlag:#3|Ort:#4|S:#5|J:#6||{%
\ifx-#1\relax\else{#1} \fi
\ifx-#2\relax\else{\bf #2}, \fi
\ifx-#4\relax\else{#4:} \fi
\ifx-#3\relax\else{#3} \fi
\ifx-#6\relax\else{#6}\fi
\ifx-#5\relax\else{, #5}\fi\Par}
\setHI{[KKLV]\ }
\sfcode`.=1000

\def\alp{\alpha} 
 
\def\bet{\beta} 
 
\def\del{\delta} 
\def\eps{\epsilon} 
\def\lam{\lambda}

\def\wlam{w_\lam}

\def\lbar{{\overline\lambda}}

\def\mbar{{\overline\mu}}

\def \ltil{\widetilde \lam}
\def \mtil{\widetilde \mu}

\fontdef Ss:cmss10,,.
\font\BF=cmbx10 scaled \magstep1
\font\CSC=cmcsc10 %scaled \magstephalf

\baselineskip15pt

{\baselineskip1.5\baselineskip\rightskip0pt plus 5truecm
\leavevmode\vskip0truecm\noindent
\BF A new formula for weight multiplicities and characters}
\vskip1truecm
\leftline{{\CSC Siddhartha Sahi}%
\footnote*{\rm This work was supported by an NSF grant.}}
\leftline{Department of Mathematics, Rutgers University, 
New Brunswick NJ 08903, USA
}
\bigskip
\beginsection intro. Introduction

The weight multiplicities of a representation of a 
simple Lie algebra $\fg$ are the dimensions of eigenspaces 
with respect to a Cartan subalgebra $\fh$. In this paper we 
give a new formula for these multiplicities.

Our formula expresses the multiplicities as sums of positive
rational numbers. Thus it is very different from the
classical formulas of Freudenthal \cite{F} and Kostant 
\cite{Ks}, which express them as sums of positive and 
negative integers. It is also quite different from recent 
formulas due to Lusztig \cite{L1} and Littelmann \cite{Li}. 

For example, for the multiplicity of the next-to-highest
weight in the $n$-dimensional representation of $\fs\fl_2$, 
we get the following expression (which sums to 1):
$${1\over (1)(2)}+{1\over (2)(3)}+\cdots+
{1\over (n-1)(n)}+ {1\over n }$$ 
 
The key role in our formula is played by the {\it dual\/} affine 
Weyl group. 

Let $V_0,(,)$ be the real Euclidean space spanned by
the root system $R_0$ of $\fg$ and let $V$ be the space of 
affine linear functions on $V_0$.  We shall identify $V$ with 
$\RR\delta\oplus V_0$ via the pairing
$(r\delta+x,y)=r+(x,y) \quad \|for|\; r\in \RR, x,y\in V_0.$

The dual affine root system is $R=\{m\delta+\alp^\vee\mid m\in\ZZ,
\alp \in R_0\}\subseteq V$ where $\alp^\vee$ means $
2\alp\over(\alp,\alp)$
as usual. Fix a positive subsystem $R_0^+\subseteq R_0$ with base 
$\{\alp_1,\cdots,\alp_n\}$  and let $\bet$ be the highest 
{\it short\/} root. Then a base for $R$ is given by 
$a_0=\del-\bet^\vee,a_1=\alp_1^\vee, \cdots, a_n=\alp_n^\vee$, 
and we write $s_i$ for the (affine) reflection 
about the hyperplane $\{x\mid (a_i,x)=0\}\subseteq V_0$.

The dual affine Weyl group is the Coxeter group $W$ generated
by $s_0,\cdots, s_n$, and the finite Weyl group is the 
subgroup $W_0$ generated by $s_1,\cdots,s_n.$ For $w\in W$, 
its {\it length\/}  is the length of a reduced (i.e. shortest) 
expression of $w$ in terms of the $s_i$. The group
$W$ acts on the weight lattice $P$ of $\fg$, and each orbit
contains a unique (minuscule) weight from the set:
$$\cO:= \{\lam\in P\mid (\alp^\vee,\lam)= 0\;\|or| 
\;1\; \forall \alp \in R^+\}.$$ 

\Definition: For each $\lam$ in $P$, we define 
\item{(1)} $\ltil:= \lam+{1\over2}
\sum_{\alp\in R_0^+}\eps_{(\alp^\vee,\lam)}\alp$ 
where, for $t\in\RR$, $\eps_t$ is $1$ if $t>0$ and $-1$ if $t\le0$
\item{(2)} $\wlam:=$  unique shortest element in $W$ such that 
$\lbar:= \wlam\cdot\lam \in \cO$

We fix a reduced expression $s_{i_1}\cdots s_{i_m}$ for $w_\lam$, 
and for each $J\subseteq \{1,\cdots,m\}$ we define
\item{(3)} $w_J:=$ the element of $W$ obtained by {\it deleting\/} 
$s_{i_j},j\in J$ from the product $s_{i_1}\ldots s_{i_m}$
\item{(4)} $c_J:= \prod_{j\in J}c_j$ where $c_j:=
(a_{i_j},\widetilde{\lam_{(j)}})^{-1}$ and
$\lam_{(j)}:= s_{i_{j-1}}\ldots s_{i_1}\cdot \lbar$

%We can now describe our formula for weight multiplicities.
Let $P^+\subset P$ be the cone of dominant weights;
and for $\lam \in P^+$, let $V_\lam$ be the irreducible 
representation of $\fg$ with 
highest weight $\lam$. 

\Theorem main. For $\lam$ in $P^+$ and $\mu$ in $P$, the 
multiplicity $m_\lam(\mu)$ of $\mu$ in $V_\lam$ is given by 
$m_\lam(\mu):={|W_0\cdot \lam|\over |W_0\cdot \mu|}
\sum_Jc_J$ where the summation is over all $J$ such that 
$w_J^{-1}\cdot \lbar$ is in $W_0\cdot\mu$. 

(We shall prove in \cite{cJpos} below that the $c_J$'s are positive.)

For $\mu$ in $P$, let $e^\mu$ denote the function $x\mapsto 
e^{(\mu,x)}$ on $V_0$. Then $W$ acts on the $e^\mu$'s by virtue of
its action on $P$, i.e. $s_i e^\mu= e^{s_i\cdot \mu}$,
and \cite{main} is equivalent to the following formula for
the character $\chi_\lam:=\sum_{\mu} m_\lam(\mu)e^\mu$ of 
$V_\lam$:

\Theorem charac. We have $\displaystyle \chi_\lam=
{|W_0\cdot \lam|\over |W_0|}\sum_{w\in W_0} w
(s_{i_m}+c_m) \cdots (s_{i_1}+c_1)e^{\lbar}$.

We will obtain \cite{charac} as a consequence of a more general
result; namely an analogous formula for the generalized Jacobi
polynomials $P_\lam$ of Heckman and Opdam. 

For the definition and properties of $P_\lam$ we refer the reader 
to \cite{HS} and \cite{Op}. We recall here that the $P_\lam$ depend 
on certain parameters $k_\alp,\alp\in R_0$ such that $k_{w\cdot 
\alp}=k_\alp$ for all $w\in W_0$. For special values of
$k_\alp$, the $P_\lam$ can be interpreted as spherical functions on
a compact symmetric space. In particular, in the limit as all
$k_\alp \rightarrow1$ we have $P_\lam\rightarrow \chi_\lam$.

\Definition: In the context of the previous definition, for 
$\lam$ in $P$ we {\it redefine\/} 
\item{($1'$)} $\ltil:= \lam+{1\over2}
\sum_{\alp\in R_0^+}k_\alp \eps_{(\alp^\vee,\lam)}\alp$
\item{($4'$)} $c_j= k_{i_j}(a_{i_j},\widetilde{\lam_{(j)}})^{-1}$
where $k_0=k_\bet$ and $k_i=k_{\alp_i}$ for $i\ge 1$.

\Theorem Plam. For $\lam$ in $P^+$, and $c_j$ as above, the 
Heckman-Opdam polynomial $P_\lam$ is given by the same formula 
as in \cite{charac}.

For $\lam$ in $P^+$ define $c_\lam:= 
{|W_0|\over |W_0\cdot \lam|}\prod_j (a_{i_j},\widetilde{\lam_{(j)}})$,
and let $\cP:=\ZZ_+[k_\alp]$ be the set of polynomials in the 
parameters $k_\alp$ with non-negative integral coefficients. 
Then we prove

\Theorem posit. $c_\lam$ is in $\cP$, as are all coefficients of
$c_\lam P_\lam$.

\cite{posit} is a generalization of the main result of
\cite{KS} to arbitrary root systems.

Our proof depends on three fundamental ideas in
the ``new'' theory of special functions.

The first idea, due to Macdonald, Heckman, Opdam, and others, 
is that one can treat root multiplicities on a symmetric space 
as parameters. 

The second idea, due to Dunkl and Cherednik, is that radial
parts of invariant differential operators on symmetric spaces 
can be written as polynomials in certain commuting first order 
differential-reflection operators --- the Cherednik operators. 

The third idea is the method of intertwiners for Cherednik 
operators. This was developed in \cite{KS}, \cite{Kn}, 
\cite{S1}, and in \cite{C2}, and can be regarded as the double
affine analog of Lusztig's fundamental relation \cite{L2}
in the affine Hecke algebra.

Using the intertwiners of \cite{C2} and \cite{S2}, our
results can be extended to the context of Macdonald
polynomials and to the six parameter Koornwinder polynomials.
These intertwiners correspond to the affine Weyl group (rather 
than the dual affine Weyl group) and hence are {\it not\/} 
appropriate for the present context. We shall discuss them
elsewhere \cite{S3}.

\beginsection prelim. Preliminaries

The results of this section are due to Cherednik [C1], 
Heckman, and Opdam [O]. 

Let $\FF=\RR(k_\alp)$ be the field of 
rational functions in the parameters $k_\alp$, and let 
$\cR$ be the $\FF$-span of $\{e^\lam\mid\lam\in P\}$ regarded 
as a $W$-module.

\Definition: For $y\in V_0$, the Cherednik operator $D_y$
is defined by 
$$D_y=\partial_y +\sum_{\alp\in R_0^+}(y,\alp) k_\alp
{1 \over 1-e^{-\alp}}(1-s_\alp) -(y,\rho);\quad 
\|where|\; 
\rho:={1\over2}\sum_{\alp\in R_0^+}k_\alp\alp$$ 

Here are some basic facts about Cherednik operators from section 2 
of [O]:

\Proposition opdam.
\item{(1)} The operators $D_y$ act on $\cR$ and commute pairwise. 
\item{(2)} For $i=1,\cdots,n$, we have $s_iD_y-D_{s_i y} s_i
=-k_i(y,\alp_i)$.
\item{(3)} There is a basis $\{E_\lam\mid \lam \in P\}$ 
of $\cR$, characterized uniquely as follows:
\item{}\quad {a)} the coefficient of $e^\lam$ in $E_\lam$ is $1$,
\item{}\quad {b)} $D_y E_\lam=(y,\ltil)E_\lam$ where $\ltil$
is as in Definition $(4')$ of the introduction.
\item{(4)} For $\lam$ in $P^+$, the Heckman-Opdam polynomial 
$P_\lam$ equals ${|W_0\cdot \lam|\over |W_0|}\sum_{w\in W_0} 
w E_\lam$.
\item{(5)} For $i=1,\cdots,n$, if $s_i\cdot \lam \ne \lam$ then 
$\widetilde{s_i\cdot\lam}=s_i\cdot \ltil$. \qed

\beginsection affine. The affine reflection

In this section we prove some basic properties of the
affine reflection $s_0$.  

\Lemma short. If $\alp$ is a positive root different from $\bet$,
then $(\alp^\vee,\bet)$ equals $0$ or $1$. 

\Proof: Since $\bet$ is in $P^+$, $(\alp^\vee,\bet)$ is a nonnegative 
integer. Also, since $\bet$ is a short root we have $(\alp,\alp)\ge
(\bet,\bet)$. So by the Cauchy-Schwartz inequality we get
$$(\alp^\vee,\bet)=2{(\alp,\bet)\over (\alp,\alp)}\le 2{(\alp,\bet)\over (\alp,\alp)^{1/2}(\bet,\bet)^{1/2}}
\le2.$$ 
If $\alp\ne\bet$, then $\alp$ is not proportional to $\bet$
and the last inequality is strict. \qed

For $i=0,1,2,$ define $R^i_0=\{\alp\in R^+_0\mid 
(\alp^\vee,\bet)= i\}$, and for $\alp$ in $R_0^+$ put 
$$\alp'=\cases{ s_\bet\cdot \alp & if $\alp\in R_0^0$ \cr
 -s_\bet\cdot \alp & if $\alp \in R_0^1\cup R_0^2$ }$$

\Lemma three. The involution $\alp\mapsto \alp'$ acts trivially
on $R_0^0$ and $R_0^2$, and permutes $R_0^1$. 

\Proof: For $\alp$ in $R_0^1$, we have
$(\alp'^\vee,\bet)=(\alp^\vee,-s_\bet\cdot\bet)=(\alp^\vee,\bet)=1$,
which implies that $\alp'$ is a (positive) root in $R_0^1$. 
The assertions about $R_0^0$ and $R_0^2=\{\bet\}$ are obvious. \qed

\Lemma s0ltil. For $\lam$ in $P$, if $s_0\cdot \lam \ne \lam$ 
then $\widetilde{s_0\cdot\lam}=s_0\cdot \ltil$. 

\Proof: We compute $s_0\cdot\ltil=\bet+s_\bet\ltil$ using 
\cite{three} and $k_{\alp}=k_{\alp'}$. This gives
$$s_0\cdot\ltil= \bet +s_\bet\cdot\lam+{1\over2}
\sum_{\alp\in R_0^0}k_\alp \eps_{(\alp^\vee,\lam)}\alp
-{1\over2}\sum_{\alp\in R_0^1\cup R_0^2}k_\alp 
\eps_{(\alp'^\vee,\lam)} \alp
$$

Comparing this to the expression for $\widetilde\mu$ with
$\mu =s_0\cdot\lam$, it suffices to show that 
$$\eps_{(\alp^\vee,\mu)}=\cases{
\eps_{(\alp^\vee,\lam)} & if $\alp\in R_0^0$\cr
-\eps_{(\alp'^\vee,\lam)} & if $\alp\in R_0^1\cup R_0^2$} $$

For $\alp$ in $R_0^0$, we easily compute that 
$(\alp^\vee,\mu)= (\alp^\vee,\lam)$.

For $\alp$ in $R_0^1$ we get 
$(\alp^\vee,\mu)=(\alp^\vee,\bet+s_\bet \cdot \lam)=
1-(\alp'^\vee,\lam).$
Being an integer $(\alp'^\vee,\lam)$ is either $\le0$
or $\ge1$. In either case
we get $\eps_{(\alp^\vee,\mu)}=-\eps_{(\alp'^\vee,\lam)}$.

Finally, for $\alp$ in $R_0^2$ we have $\alp=\alp'=\bet$
and $(\bet^\vee,\mu) = 2-(\bet^\vee,\lam)$. 
Now $s_0\lam\ne\lam$ implies
that $(\bet^\vee,\lam)\ne 1$, thus we have either
$(\bet^\vee,\lam)\ge2$ or $(\bet^\vee,\lam)\le0$.
In either case we get $\eps_{(\bet^\vee,\lam)}=
\eps_{(\bet^\vee,\lam)}=
-\eps_{(\bet^\vee,\mu)}.$ \qed
\beginsection intertwine. The intertwining relation

Dualizing the action $y\mapsto w\cdot y$ of $W$ on $V_0$, 
we get a representation $v\mapsto wv$ of $W$ on $V$ satisfying
$(wv, y)=(v, w^{-1}\cdot y).$ For $y$ in $V_0$ and $w$ in $W_0$, 
we have $wy=w\cdot y$. The affine reflection $s_0$ acts 
on $V$ by $$s_0(r\del+y) = (y,\bet)\del+ r \del+s_\bet y. $$ 

For $v=r\del+y$ in $V$ we define the {\it affine\/} Cherednik 
operator simply by putting $D_v=D_y+rI$ where $I$ is the identity 
operator. 
From (2) of \cite{opdam} we know the intertwining 
relations between the (affine) Cherednik operators and 
$s_1,\cdots,s_n$. In this section we prove the following 
intertwining relation between these operators and $s_0$:

\Proposition inter. For $v=r\delta+y$ in $V$, we have 
$D_vs_0- s_0D_{s_0 v}= k_\bet(y,\bet)$

\Proof: Let us write $N_\alp$ for ${1\over 1-e^{-\alp}}(1-s_\alp)$, so that 
$$D_v= \partial_y+\sum k_\alp (y,\alp) N_\alp -(y,\rho)+r$$

Since $s_\bet N_\alp= N_{s_\bet\cdot \alp}s_\bet$, and $s_\bet\partial_y= 
\partial_{s_\bet y}s_\bet$, we get  
$$s_\bet D_v s_\bet = \partial_{s_\bet y}+\sum_{\alp\in R_0^+} 
k_\alp (y,\alp) N_{s_\bet\cdot\alp} -(y,\rho)+r$$

Now partitioning $R_0^+=R_0^0\cup R_0^1\cup R_0^2$ and using \cite{three} 
we get:
$$s_\bet D_vs_\bet= \partial_{s_\bet y}+
\sum_{\alp\in R_0^0} k_\alp (s_\bet y,\alp) N_{\alp} 
-\sum_{\alp\in R_0^1\cup R_0^2} k_\alp (s_\bet y,\alp) N_{-\alp} 
-(y,\rho)+r$$

The following identities are easy to check:
\item{1)} $e^\bet\partial_{s_\bet y}e^{-\bet} = 
\partial_{s_\bet y}+ (y,\bet)$ 
\item{2)} $e^\bet N_\alp e^{-\bet}=N_\alp $ for $\alp \in R_0^0$
\item{3)} $e^\bet N_{-\alp}e^{-\bet} =1-N_\alp$ for $\alp \in R_0^1$
\item{4)} $e^\bet N_{-\bet}e^{-\bet} =1-N_\bet+s_0$ 

Using these we get the following formula for $s_0 D_vs_0=
e^\bet(s_\bet D_v s_\bet )e^{-\bet}$:
$$%\eqalign{ &=
\partial_{s_\bet y} 
+(y,\bet)+ \sum_{\alp\in R_0^+} k_\alp (s_\bet y,\alp) N_\alp 
-\sum_{\alp\in R_0^1\cup R_0^2} k_\alp (s_\bet y,\alp) 
%\cr &\quad\quad 
-k_\bet(s_\bet y,\bet)s_0-(y,\rho)+r$$

Since $\sum_{\alp\in R_0^1\cup R_0^2} k_\alp(s_\bet y,\alp) = 
( s_\bet y,\rho-s_\bet\cdot\rho)= (s_\bet y,\rho)-( y,\rho)$, we get 
$$s_0 D_vs_0 =D_{s_\bet y} +(y,\bet) -k_\bet(s_\bet y,\bet)s_0+r
=D_{s_0 v}+k_\bet(y,\bet)s_0$$

The result follows. \qed

\beginsection result. The Heckman-Opdam polynomials

Let $E_\lam$ be as in \cite{opdam}.

\Proposition jacobi. The polynomials $E_\lam$ satisfy the following
recursions:
\item{1)} $E_\lam=e^\lam$ for $\lam \in \cO$
\item{2)} If $s_i\cdot\lam\ne\lam$ then 
$\left(s_i+ {k_i\over(a_i,\ltil)}\right) E_\lam$ 
is a multiple of $E_{s_i\cdot\lam}$

\Proof: For 1) we check simply that $D_y e^\lam=(y,\ltil)e^\lam$,  
 using the identity 
$$N_\alp e^\lam=\cases{e^\lam & if $(\alp^\vee,\lam)=1$\cr 
0 & if $(\alp^\vee,\lam)=0$}.$$ 

For 2), we write $F$ for $\left(s_i+ {k_i\over(a_i,\ltil)}\right)
 E_\lam$ and first consider $i\ne0$. Then for $y$ in $V_0$,
using (2) of \cite{opdam}, we get
$$ D_yF =\left(s_iD_{s_i y}-k_i(y,\alp_i) +
{k_i\over(a_i,\ltil)}D_y \right)E_\lam
=\left( (s_i y,\ltil)s_i
 +k_i {(y,\ltil)\over (a_i,\ltil)} -k_i(y,\alp_i)\right)E_\lam
$$
Since $(y,\ltil)-(y,\alp_i)(a_i,\ltil) =
 (s_iy,\ltil)$, using (5) of \cite{opdam} we get
$$D_yF =(s_iy,\ltil)F=(y,s_i\cdot \ltil)F=
(y,\widetilde{s_i\cdot\lam})F.$$
This proves 2) for $i\ne0$. For $i=0$ we use the previous 
proposition to get
$$ D_yF =\left(s_0D_{s_0 y}+k_\bet(y,\bet) +
{k_\bet\over(a_0,\ltil)}D_y \right)E_\lam
=\left( (s_0 y,\ltil)s_0 +k_\bet {(y,\ltil)\over (a_0,\ltil)} 
+k_\bet(y,\bet)\right)E_\lam$$
This time using $(y,\ltil)+ (a_0,\ltil)(y,\bet)=(s_0y,\ltil)$, 
and \cite{s0ltil} we get
$$D_yF =(s_0y,\ltil)F=(y,s_0\cdot \ltil)F=
(y,\widetilde{s_0\cdot \lam})F.$$
This completes the proof of 2) for $i=0$. \qed

\Corollary Elam. For $\lam$ in $P$, and $c_i$ as in 
Definition $(4')$ of the introduction, we have 
$$E_\lam=(s_{i_m}+c_m) \cdots (s_{i_1}+c_1) e^{\lbar}.$$

\Proof: By the minimality of $w_\lam$, if $w$ is a {\it proper\/} 
subexpression of $w_\lam^{-1}= s_{i_m}\cdots s_{i_1}$ then 
$w\cdot \lbar\ne\lam$. This means that the coefficient of $e^\lam$
in $(s_{i_m}+c_m)\cdots (s_{i_1}+c_1) e^{\lbar}$ is $1$.
The result now follows from the previous proposition. \qed

\Proof: (of \cite{Plam}) This follows from
\cite{Elam} and \cite{opdam} (4). \qed

\beginsection positivity. Positivity

Let $\cP_1\subset \cP$ be the set of polynomials of 
degree $\le1$, with non-negative integral coefficients,
and {\it positive\/} constant term. 

For $\lam$ in $P$, let $a_{i_j}$ and $\widetilde {\lam_{(j)}}$ 
be as in Definition $(4')$ of the introduction.

\Proposition pos. For each $j=1,\cdots,m$,
 $(a_{i_j},\widetilde {\lam_{(j)}})$ belongs to $\cP_1$. 

\Proof: Fix $j$ and write $\mu=\lam_{(j)}$, $i=i_j$, and 
$w=s_{i_1}\cdots s_{i_{j-1}}$. We need to show that $(a_i,\mtil)$ 
has positive constant term and non-negative integral coefficients.

The lengths of $w$ and $ws_i$ must be $j-1$ and $j$ respectively,
since otherwise we could shorten the expression $s_{i_1}\cdots 
s_{i_m}$ for $w_\lam$. By a standard argument (\cite{Hu}, Ch. 5) 
this implies that $w(a_i)$ is a positive (affine) coroot in $R^+$. 
Since $\lbar=\mbar$ is minuscule, we conclude
$$0\le(w(a_i),\mbar)=(a_i,w^{-1}\cdot\mbar)=(a_i,\mu)$$
If $(a_i,\mu)$ were $0$ then $\lam_{(j+1)}=s_i\cdot \mu=
\mu=\lam_{(j)}$ and we could shorten the expression for 
$w_\lam$ by dropping $s_{i_j}$. 
This shows that $(a_i,\mu)$, which is the constant term of 
$(a_i,\mtil)$, is positive.

If $i=0$, the non-constant part of $(a_0,\mtil)$ is
$$\displaystyle -{1\over2}\sum_{\alp\in R_0^+}k_\alp
\eps_{(\alp^\vee,\mu)}(\bet^\vee,\alp),$$ 
and we consider separately the contributions of $R_0^0,R_0^1,
R_0^2$. 
For $\alp$ in $R_0^0$ the contribution is $0$. 
For $\alp=\bet$ in $R_0^2$, we get the term $-\eps_{(\bet^\vee,\mu)}k_\bet$. 
By the first part, $(a_0,\mu)$ is a positive integer. Hence
$(\bet^\vee,\mu)=1-(a_0,\mu)\le0$, which implies 
that $-\eps_{(\bet^\vee,\mu)}=1$.
The roots in $R_0^1$ can be grouped in pairs 
$\{\alp,-s_\bet\cdot\alp\}$ and the contribution of
such a pair is
$$-k_\alp{\eps_{(\alp^\vee,\mu)}+\eps_{(-s_\bet \alp^\vee,\mu)}\over2}
(\bet^\vee,\alp).$$ 
Now $(\bet^\vee,\alp)$ is positive, so the coefficient above
is a non-negative integer, unless
$(\alp^\vee,\mu)$  and $(-s_\bet \alp^\vee,\mu)$ are both $>0$. 
But in this case, we would get
$$0<(\alp^\vee,\mu)-(s_\bet \alp^\vee,\mu)
=(\alp^\vee,\mu-s_\bet\cdot \mu)= (\bet^\vee,\mu)
(\alp^\vee,\bet)\le0$$ which is a contradiction.

The argument is similar if $i>0$. The non-constant part of $(a_i,\mtil)$ is 
$${1\over2}\sum_{\alp\in R_0^+}k_\alp
\eps_{(\alp^\vee,\mu)}(a_i,\alp).$$ 
To compute this we divide $R_0^+$ into three disjoint sets 
consisting of $\{\alp_i\}$, $\{$the roots orthogonal to 
$\alp_i\}$, and $\{$the remaining positive roots$\}$. 
For $\alp=\alp_i$, we get the coefficient $\eps_{(a_i,\mu)}$ 
which is $1$ since $(a_i,\mu)>0$ by the first part.
If $\alp$ is orthogonal to $\alp_i$ then the coefficient 
is zero. Finally, the remaining positive roots can be grouped 
into pairs $\{\alp, s_i\cdot \alp\}$ where we may assume 
that $(\alp^\vee,\alp_i)>0$. The contribution of 
each such pair is
$$k_\alp{\eps_{(\alp^\vee,\mu)}-\eps_{(s_i \alp^\vee,\mu)}\over2}
(a_i,\alp).$$ 
Now $(\alp^\vee,\alp_i)>0$ implies $(a_i,\alp)>0$. Therefore
this coefficient is a non-negative integer, unless
$(\alp^\vee,\mu) \le 0$ and $(s_i \alp^\vee,\mu)>0$. But if
this were the case then we would have
$$0> (\alp^\vee,\mu)-(s_i\alp^\vee,\mu)= 
(\alp^\vee,\mu-s_i\cdot\mu)=(a_i,\mu)(\alp^\vee,\alp_i)
>0,$$
which is a contradiction. \qed

\Proof: (of \cite{posit}) This follows from
\cite{Plam} and \cite{pos}. \qed

Setting all the $k_\alp$'s equal to $1$ in \cite{pos}, we deduce

\Corollary cJpos. The constants $c_j$ and $c_J$ in
Theorems \ncite{main} and \ncite{charac} are positive. 
\qed

\beginsection References. References

\baselineskip12pt
\parskip2.5pt plus 1pt
\hyphenation{Hei-del-berg}
\def\L|Abk:#1|Sig:#2|Au:#3|Tit:#4|Zs:#5|Bd:#6|S:#7|J:#8||{%
\edef\TEST{[#2]}
\expandafter\ifx\csname#1\endcsname\TEST\relax\else
\immediate\write16{#1 had already been defined!}\fi
\expandwrite\AUX{\neverexpand\ref{#1}{\TEST}}
\HI{[#2]}
\ifx-#3\relax\else{#3}: \fi
\ifx-#4\relax\else{#4}{\sfcode`.=3000.} \fi
\ifx-#5\relax\else{\it #5\/} \fi
\ifx-#6\relax\else{\bf #6} \fi
\ifx-#8\relax\else({#8})\fi
\ifx-#7\relax\else, {#7}\fi\Par}

\def\B|Abk:#1|Sig:#2|Au:#3|Tit:#4|Reihe:#5|Verlag:#6|Ort:#7|J:#8||{%
\edef\TEST{[#2]}
\expandafter\ifx\csname#1\endcsname\TEST\relax\else
\immediate\write16{#1 hat sich geaendert!}\fi
\expandwrite\AUX{\neverexpand\ref{#1}{\TEST}}
\HI{[#2]}
\ifx-#3\relax\else{#3}: \fi
\ifx-#4\relax\else{#4}{\sfcode`.=3000.} \fi
\ifx-#5\relax\else{(#5)} \fi
\ifx-#7\relax\else{#7:} \fi
\ifx-#6\relax\else{#6}\fi
\ifx-#8\relax\else{ #8}\fi\Par}

\def\Pr|Abk:#1|Sig:#2|Au:#3|Artikel:#4|Titel:#5|Hgr:#6|Reihe:{%
\edef\TEST{[#2]}
\expandafter\ifx\csname#1\endcsname\TEST\relax\else
\immediate\write16{#1 hat sich geaendert!}\fi
\expandwrite\AUX{\neverexpand\ref{#1}{\TEST}}
\HI{[#2]}
\ifx-#3\relax\else{#3}: \fi
\ifx-#4\relax\else{#4}{\sfcode`.=3000.} \fi
\ifx-#5\relax\else{In: \it #5}. \fi
\ifx-#6\relax\else{(#6)} \fi\PrII}
\def\PrII#1|Bd:#2|Verlag:#3|Ort:#4|S:#5|J:#6||{%
\ifx-#1\relax\else{#1} \fi
\ifx-#2\relax\else{\bf #2}, \fi
\ifx-#4\relax\else{#4:} \fi
\ifx-#3\relax\else{#3} \fi
\ifx-#6\relax\else{#6}\fi
\ifx-#5\relax\else{, #5}\fi\Par}
\setHI{[ABC]\ }
\sfcode`.=1000

\L|Abk:C1|Sig:C1|Au:Cherednik, I.|Tit:Double affine Hecke algebras
and Macdonald's conjectures|Zs:Ann. Math|Bd:141|S:191--216|J:1997||

\L|Abk:C2|Sig:C2|Au:Cherednik, I.|Tit:Intertwining operators of double
affine Hecke algebras|Zs:Selecta Math|Bd:-|S:-|J:1997||

\L|Abk:F|Sig:F|Au:Freudenthal, H.|Tit:Zur Berechnung der
chararktere der halbeinfachen Lieschen gruppen|Zs:Indag. 
Math.|Bd:16|S:369--376|J:1954||

\B|Abk:HS|Sig:HS|Au:Heckman, G.; Schlichtkrull, H.|Tit:Harmonic 
Analysis and Special Functions on Symmetric 
Spaces|Reihe:-|Verlag:Academic Press.|Ort:-|J:1994||

\B|Abk:Hu|Sig:Hu|Au:Humphreys, J.|Tit:Reflection Groups
and Coxeter Groups|Reihe:-|Verlag:Cambridge Univ. 
Press.|Ort:-|J:1990||

\L|Abk:Kn|Sig:Kn|Au:Knop, F.|Tit:Integrality of two variable
Kostka functions|Zs:J. Reine Ang. Math.|Bd:482|S:177--189|J:1997||

\L|Abk:KS|Sig:KS|Au:Knop, F.; Sahi, S.|Tit:A recursion 
and a combinatorial formula for Jack polynomials|Zs:Invent. 
math.|Bd:128|S:9--22|J:1997||

\L|Abk:Ko|Sig:Ko|Au:Koornwinder, T.|Tit:Askey-Wilson polynomials
for root systems of type $BC$|Zs:Contemp.~Math.|Bd:138|S:189-204|J:1992||

\L|Abk:Ks|Sig:Ks|Au:Kostant, B.|Tit:A formula for the multiplicity 
of a weight|Zs:Trans. Amer. Math. Soc.|Bd:93|S:53--79|J:1959||

\L|Abk:Li|Sig:Li|Au:Littelmann, P.|Tit:Paths and root operators
in representation theory|Zs:Ann. of Math.|Bd:142|S:499--525|J:1995||

\L|Abk:L1|Sig:L1|Au:Lusztig, G.|Tit:Singularities, character 
formulas and a $q$-analog of weight 
multiplicities|Zs:Ast\'erisque|Bd:101-102|S:208--299|J:1983||

\L|Abk:L2|Sig:L2|Au:Lusztig, G.|Tit:Affine Hecke algebras and
their graded version|Zs:J. Amer. Math. Soc.|Bd:2|S:599--635|J:1989||

\L|Abk:M|Sig:M|Au:Macdonald, I.|Tit:Affine Hecke algebras and 
orthogonal polynomials|Zs:S\'em. Bourbaki|Bd:797|S:1--18|J:1995||

\L|Abk:Op|Sig:O|Au:Opdam, E.|Tit:Harmonic analysis 
for certain representations of graded Hecke 
algebras|Zs:Acta Math.|Bd:175|S:75--121|J:1995||

\L|Abk:S1|Sig:S1|Au:Sahi, S.|Tit:Interpolation, 
integrality, and a generalization of Macdonald's 
polynomials|Zs:IMRN|Bd:10|S:457--471|J:1996||

\L|Abk:S2|Sig:S2|Au:Sahi, S.|Tit:Nonsymmetric Koornwinder 
polynomials and duality|Zs:-|Bd:-|S:-|J:preprint||

\L|Abk:S3|Sig:S3|Au:Sahi, S.|Tit:The Bruhat order,
Macdonald polynomials, and positivity for Heckman-Opdam
polynomials|Zs:-|Bd:-|S:-|J:preprint||
\bye